\documentclass[12pt]{article} 
\usepackage{amsfonts,amsmath,latexsym,amssymb,mathrsfs,amsthm,comment}
\usepackage{slashbox}
\usepackage{caption}

\usepackage{pgfplots}

\let\OLDthebibliography\thebibliography
\renewcommand\thebibliography[1]{
  \OLDthebibliography{#1}
  \setlength{\parskip}{3pt}
  \setlength{\itemsep}{0pt plus 0.3ex}
}

\def\numberlikeadb{\global\def\theequation{\thesection.\arabic{equation}}}
\numberlikeadb
\newtheorem{theorem}{Theorem}[section]

\newtheorem{conjecture}[theorem]{Conjecture}

\usepackage{lscape}
\usepackage{caption}
\usepackage{multirow}

\usepackage{color}
\newcommand{\gr}[1]{{\color{black} #1}}
\begin{document}

\title{The basic distributional theory for the product of zero mean correlated normal random variables}
\author{Robert E. Gaunt\footnote{Department of Mathematics, The University of Manchester, Oxford Road, Manchester M13 9PL, UK, robert.gaunt@manchester.ac.uk}}

\date{} 
\maketitle

\vspace{-5mm}

\begin{abstract}The product of two zero mean correlated normal random variables\gr{, and more generally the sum of independent copies of such random variables,} has received much attention in the statistics literature and appears in many application areas. However, many important distributional properties are yet to be {recorded}. This \gr{review} paper fills this gap by providing the basic distributional theory for the \gr{sum of independent copies of the} product of two zero mean correlated normal random variables. \gr{Properties covered include probability and cumulative distribution functions, generating functions, moments and cumulants, mode and median, Stein characterisations, representations in terms of other random variables, and a list of related distributions. We also review how the product of two zero mean correlated normal random variables arises naturally as a limiting distribution, with an example given for the distributional approximation of double Wiener-It\^{o} integrals.}
\end{abstract}


\noindent{{{\bf{AMS 2010 Subject Classification:}}} Primary 60E05; 62E15

\section{Introduction}

Let $(X,Y)$ be a bivariate normal random vector with zero mean vector, variances $(\sigma_X^2,\sigma_Y^2)$ and correlation coefficient $\rho$. The distribution of the product $Z=XY$ arises in numerous applications, with recent examples including statistics of Lagrangian power in two-dimensional turbulence (Bandi and Connaughton \cite{bandi}), product confidence limits for indirect effects (MacKinnon et al.\ \cite{mac2}) and statistical mediation analysis (MacKinnon \cite{mac1}). The mean $\overline{Z}_n=n^{-1}(Z_1+Z_2+\cdots+Z_n)$, where $Z_1, Z_2,\ldots, Z_n$ are independent and identical copies of $Z$, has also found application in areas such as electrical engineering (Ware and Lad \cite{ware}), astrophysics (Mangilli, Plaszczynski and Tristram \cite{mpt15}, Watts et al.\ \cite{watts}) and quantum cosmology (Grishchuk \cite{g96}).


The exact distribution of the product $Z=XY$ has been studied since 1936 (Craig \cite{craig}), with subsequent contributions coming from Aroian \cite{aa1}, Aroian, Taneja and Cornwell \cite{aa2}, Bandi and Connaughton \cite{bandi}, Haldane \cite{hh1}, Meeker et al$.$ \cite{mm1}; see Nadarajah and Pog\'{a}ny \cite{np16} for an overview of these and further contributions. Recently, Nadarajah and Pog\'{a}ny \cite{np16} used characteristic functions to find exact formulas for the probability density function (PDF) of the product $Z$ and more generally the mean $\overline{Z}_n$; these formulas are starting to become well-known and have found many recent applications. Independently, in the physics literature, Grishchuk \cite{g96} had obtained the formula for the PDF of $Z$ and \gr{Mangilli, Plaszczynski and Tristram \cite{mpt15}} then obtained the more general formula for the PDF of $\overline{Z}_n$. Much earlier, a formula for PDF of $Z$ was given in the book of Springer \cite[equation (4.8.22)]{springer}, but there was an unfortunate typo. 

Despite the interest in the distributions of $Z$ and $\overline{Z}_n$, many important distributional properties beyond the PDF are not recorded in the literature. \gr{Recent} work of Gaunt \cite{gaunt prod} identified $Z$ and $\overline{Z}_n$ as variance-gamma (VG) random variables (leading to a simple independent derivation of their PDFs), which means that some of their distributional properties can be inferred from results for the VG distribution; see, for example, Chapter 4 of the book of Kotz, Kozubowski and Podg\'{o}rski \cite{kkp01}. However, searching through the VG literature \gr{can be difficult}; it is tedious to convert results for the VG distribution to results for the distributions of $Z$ and $\overline{Z}_n$;
and results that exploit the special structure of the random variables $Z$ and $\overline{Z}_n$ are not available. 

In this \gr{review}, we aim to fill a gap in the literature by providing the basic distributional theory for the  product of two zero mean correlated normal random variables. We present results for the more general mean $\overline{Z}_n$, with results for the product $Z$ following on setting $n=1$. The end result is that many of the most important distributional properties of $Z$ and $\overline{Z}_n$ are now collected in a single reference. \gr{Most results are already explicitly stated in the literature or can be readily deduced from the fact that $\overline{Z}_n$ has a VG distribution, in which case we provide references. Other results stated in this paper have to the best knowledge of the author not appeared in the literature, in which case we provide straightforward and concise derivations.} 

The \gr{distributional} properties covered include: formulas for the PDF (Section \ref{sec2.1}), list of related distributions (Section \ref{sec2.3}), the cumulative distribution function (Section \ref{sec2.2}), generating functions and infinite divisibility (Section \ref{sec2.4}), representations in terms of other random variables (Section \ref{sec2.5}), Stein characterisation (Section \ref{sec2.6}), moments and cumulants (Section \ref{sec2.7}), and mode and median (Section \ref{sec2.8}). \gr{This list covers some of the most basic and important properties of a probability distribution, but is not comprehensive; we do not cover multivariate extensions, connections to the Wishart distribution or inference methods, to name a few topics. In Section \ref{sec3}, we review the topic of the role of the distributions of $Z$ and $\overline{Z}_n$ as limiting distributions, in which there has been  recent interest. In particular, we present some quantitative limit theorems concerning the distributional approximation of double Wiener-It\^{o} integrals.} Basic properties of the modified Bessel function of the second kind will be needed throughout the paper, and are collected in Appendix \ref{appa}. 

\vspace{3mm}

\gr{\noindent{\emph{Notation}.} To simplify formulae, we define $s_n:=\sigma_X\sigma_Y/n$ and $s:=s_1=\sigma_X\sigma_Y$.}





\section{Distributional theory for the product of zero mean correlated normal random variables}

\subsection{Probability density function}\label{sec2.1}

For derivations of the following formulas for the PDFs of $Z$ and $\overline{Z}_n$ we refer the reader to \gr{Mangilli, Plaszczynski and Tristram \cite{mpt15}}, Nadarajah and Pog\'{a}ny \cite{np16} and Gaunt \cite{gaunt prod, gaunt prod2}. The formula for the PDF of $Z$ was also earlier derived by Grishchuk \cite{g96} and Springer \cite{springer}, with the latter reference having a typo in the formula. 
For $x\in\mathbb{R}$,
\begin{equation}\label{pdf1}f_Z(x)=\frac{1}{\pi s\sqrt{1-\rho^2}}\exp\bigg(\frac{\rho x}{s(1-\rho^2)}\bigg)K_0\bigg(\frac{|x|}{s(1-\rho^2)}\bigg), 
\end{equation}
and, for $n\geq1$,
\begin{equation}\label{pdf2}f_{\overline{Z}_n}(x)=\frac{2^{(1-n)/2}|x|^{(n-1)/2}}{s_n^{(n+1)/2}\sqrt{\pi(1-\rho^2)}\Gamma(n/2)}\exp\bigg(\frac{\rho  x}{s_n(1-\rho^2)} \bigg)K_{\frac{n-1}{2}}\bigg(\frac{ |x|}{s_n(1-\rho^2)}\bigg).
\end{equation}
Here $K_\nu(x)$ is a modified Bessel function of the second kind; see Appendix \ref{appa} for a definition and some standard properties.

\begin{figure}[h!]
\centering
\includegraphics[width=6.2in]{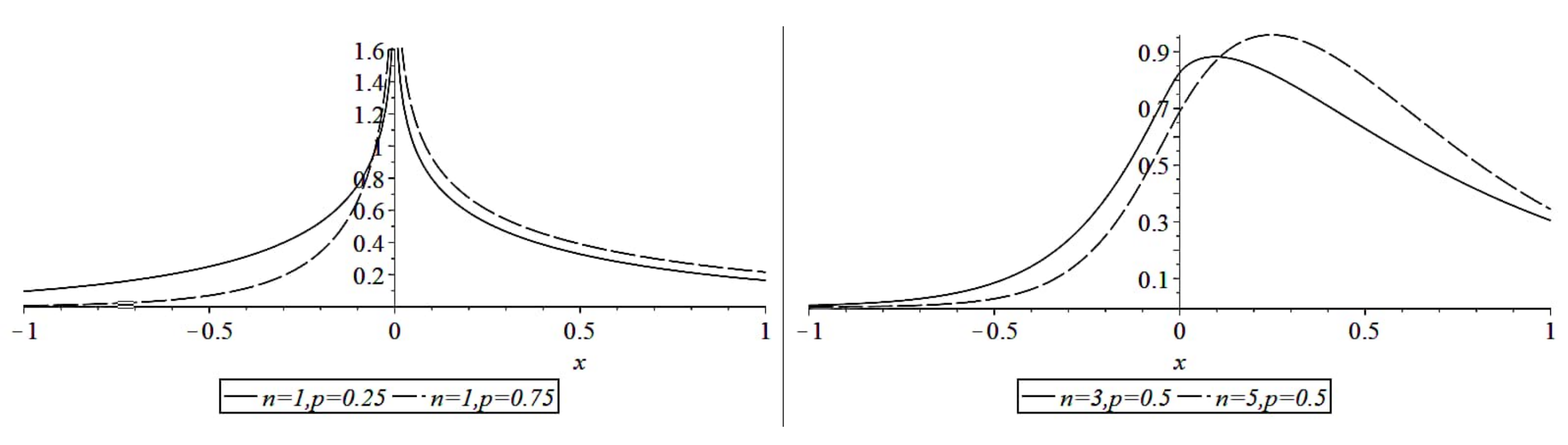}
\caption{PDFs of $Z$ and $\overline{Z}_n$ with $s=1$}
\end{figure}

The modified Bessel function in the PDFs (\ref{pdf1}) and (\ref{pdf2}) make them difficult to parse on first inspection.  We can gain some understanding from the following limiting forms.   Using the limiting form (\ref{Ktend0}), we have that 
\begin{equation}\label{pmutend}f_{\overline{Z}_n}(x)\sim\begin{cases}\displaystyle -\frac{1}{\pi s\sqrt{1-\rho^2}}\log|x|, & x\rightarrow0,\: n=1, \\
\displaystyle \frac{(1-\rho^2)^{n/2-1}}{2\sqrt{\pi}s_n}\frac{\Gamma((n-1)/2)}{\Gamma(n/2)}, &  x\rightarrow0,\:n\geq2. \end{cases}
\end{equation}
We see that the density has a singularity at the origin if $n=1$.  In fact, for all parameter values ($n\geq1$, $-1<\rho<1$, $\sigma_X,\sigma_Y>0$) the distribution of $\overline{Z}_n$ is unimodal; see Section \ref{sec2.8} for further details. For $n\geq2$, the density is bounded.  As observed by Watts et al.\ \cite{watts}, the tail behaviour of the density is obtained by applying the limiting form (\ref{Ktendinfinity}):
\begin{equation}\label{tail1}
f_{\overline{Z}_n}(x)\sim \frac{1}{(2s_n)^{n/2}\Gamma(n/2)}x^{n/2-1}\exp\bigg(-\frac{x}{s_n(1+\rho)}\bigg), \quad x\rightarrow\infty,
\end{equation}
and
\begin{equation}\label{tail2}
f_{\overline{Z}_n}(x)\sim \frac{1}{(2s_n)^{n/2}\Gamma(n/2)}(-x)^{n/2-1}\exp\bigg(\frac{x}{s_n(1-\rho)}\bigg), \quad x\rightarrow-\infty.
\end{equation}
We observe that the distribution of $\overline{Z}_n$ has semi-heavy tails.

In the case that $n$ is even, we can take advantage of a standard simplification of the modified Bessel function of the second kind (see (\ref{special})) to obtain a representation of the PDF in terms of elementary functions:
\begin{align}f_{\overline{Z}_n}(x)&=\frac{1}{(2s_n)^{n/2}\Gamma(n/2)}|x|^{n/2-1}\exp\bigg(\frac{\rho x-|x|}{s_n(1-\rho^2)}\bigg)\times\nonumber\\
\label{nsimp}&\quad\times\sum_{j=0}^{n/2-1}\frac{(n/2-1+j)!}{(n/2-1-j)!j!}\bigg(\frac{s_n(1-\rho^2)}{2|x|}\bigg)^j, \quad x\in\mathbb{R},\:n\in 2\mathbb{Z}^+.
\end{align}

\subsection{Related distributions}\label{sec2.3}

\noindent{1.} Setting $n=2$ in (\ref{nsimp}) yields
\begin{equation}\label{lapd}f_{\overline{Z}_2}(x)=\frac{1}{s}\exp\bigg(\frac{2(\rho x-|x|)}{s(1-\rho^2)}\bigg),
\end{equation}
which is the PDF of the asymmetric Laplace distribution (see Chapter 3 of Kotz, Kozubowski and Podg\'{o}rski \cite{kkp01} for a comprehensive account of its distributional theory).

\vspace{3mm}

\noindent{2.} As $\rho\uparrow1$, $\overline{Z}_n$ converges in distribution to a scaled chi-square random variable: $\overline{Z}_n\rightarrow_d s_nV$, where $V\sim \chi_{(n)}^2$, the chi-square distribution with $n$ degrees of freedom. \gr{This fact was observed by Watts et al.\ \cite{watts}.} By symmetry, $\overline{Z}_n\rightarrow_d -s_nV$, as $\rho\downarrow-1$. This is hardly surprising given the definition of $\overline{Z}_n$, and can be confirmed rigorously by letting $\rho\rightarrow\pm1$ in the formula (\ref{cfcf}) below for the characteristic function of $\overline{Z}_n$ and comparing to the chi-square characteristic function using L\'evy's continuity theorem.

\vspace{3mm}

\noindent{3.} By the central limit theorem, $\sqrt{n}(\overline{Z}_n-\rho s)\rightarrow_d N(0,s^2(1+\rho^2))$, as $n\rightarrow\infty$. Here, $\mathbb{E}[Z]=\rho s$ and $\mathrm{Var}(Z)=s^2(1+\rho^2)$ (see (\ref{mean}) and (\ref{var})).

\vspace{3mm}

\noindent{4.} The variance-gamma (VG) distribution with parameters $r > 0$, $\theta \in \mathbb{R}$, $\sigma >0$, $\mu \in \mathbb{R}$ has PDF
\begin{equation}\label{vgdef}f_{\mathrm{VG}}(x) = \frac{\mathrm{e}^{\theta (x-\mu)/\sigma^2}}{\sigma\sqrt{\pi} \Gamma(r/2)}  \bigg(\frac{|x-\mu|}{2\sqrt{\theta^2 +  \sigma^2}}\bigg)^{\frac{r-1}{2}} K_{\frac{r-1}{2}}\bigg(\frac{\sqrt{\theta^2 + \sigma^2}}{\sigma^2} |x-\mu| \bigg), \quad x\in\mathbb{R}, 
\end{equation}
If a random variable $V$ has PDF (\ref{vgdef}), we write $V\sim \mathrm{VG}(r,\theta,\sigma,\mu)$.  This parametrisation was given in Gaunt \cite{gaunt vg}.  It is similar to the parametrisation given by Finlay and Seneta \cite{finlay} and other parametrisations are given by Eberlein and Hammerstein \cite{eberlein} and Kotz, Kozubowski and Podg\'{o}rski \cite{kkp01}, with Chapter 4 of the latter reference giving the most comprehensive account of the distributional theory of the VG distribution in the literature. The VG distribution is also known as the Bessel function distribution or the McKay Type II distribution (McKay \cite{m32}), as well as the generalized Laplace distribution (Kotz, Kozubowski and Podg\'{o}rski \cite{kkp01}).


It was noted in the thesis of Gaunt \cite{gaunt thesis} that $Z\sim\mathrm{VG}(1,\rho s,s\sqrt{1-\rho^2},0)$ and later by Gaunt \cite{gaunt prod} that also 
\begin{equation}\label{vgrep}\overline{Z}_n\sim\mathrm{VG}(n,\rho s_n,s_n\sqrt{1-\rho^2},0).
\end{equation}

\vspace{3mm}

\noindent{5.} \gr{The following rather complicated exact formula for the PDF of the product of correlated normal random variables with non-zero means was derived by Cui et al.\ \cite{cui}. Let $(X,Y)$ be a bivariate normal random vector with mean vector $(\mu_X,\mu_Y)$, variances $(\sigma_X^2,\sigma_Y^2)$ and correlation coefficient $\rho$. Then, for $x\in\mathbb{R}$,
\begin{align*}f_{Z}(x)&=\exp\bigg\{-\frac{1}{2(1-\rho^2)}\bigg(\frac{\mu_X^2}{\sigma_X^2}+\frac{\mu_Y^2}{\sigma_Y^2}-\frac{2\rho(x+\mu_X\mu_Y)}{\sigma_X\sigma_Y}\bigg)\bigg\}\\
&\quad\times\sum_{n=0}^\infty\sum_{m=0}^{2n}\frac{x^{2n-m}|x|^{m-n}\sigma_X^{m-n-1}}{\pi(2n)!(1-\rho^2)^{2n+1/2}\sigma_Y^{m-n+1}}\binom{2n}{m}\bigg(\frac{\mu_X}{\sigma_X^2}-\frac{\rho \mu_Y}{\sigma_X\sigma_Y}\bigg)^m\\
&\quad\times\bigg(\frac{\mu_Y}{\sigma_Y^2}-\frac{\rho \mu_X}{\sigma_X\sigma_Y}\bigg)^{2n-m}K_{m-n}\bigg(\frac{|x|}{(1-\rho^2)\sigma_X\sigma_Y}\bigg).
\end{align*}

The following exact formula for the PDF of the product $Z_k=\prod_{j=1}^k X_j$ of independent normal random variables $X_j\sim N(0,\sigma_{X_j}^2)$, $j=1,\ldots,k$, was obtained by Springer and Thompson \cite{springer2}. Let $\sigma_k=\sigma_{X_1}\sigma_{X_2}\cdots \sigma_{X_k}$. Then, for $x\in\mathbb{R}$,
\begin{equation*}f_{Z_k}(x)=\frac{1}{(2\pi)^{k/2}\sigma_k}G_{0,k}^{k,0}\bigg(\frac{x^2}{2^k\sigma_k^2} \; \bigg| \;0,\ldots,0\bigg),
\end{equation*}
where $G_{0,k}^{k,0}(x\;|\;0,\ldots,0)$ is a Meijer $G$-function (for a definition and basic properties, see Chapter 16 of Olver et al.\ \cite{olver}).}
\gr{Now, suppose $\rho=0$, and let $\overline{Z}_{n,1},\ldots,\overline{Z}_{n,k}$ be independent copies of $\overline{Z}_{n}$. A formula for the PDF of the product $Z_{n,k}=\prod_{j=1}^k\overline{Z}_{n,j}$ can be read off from a formula of Gaunt, Mijoule and Swan \cite{gms}. For $x\in\mathbb{R}$,
\begin{equation*}f_{Z_{n,k}}(x)=\frac{n}{2^k\pi^{k/2}s^k(\Gamma(n/2))^k}G_{0,2k}^{2k,0}\bigg(\frac{n^2x^2}{2^{2k}s^k} \; \bigg| \;\frac{n-1}{2},\ldots,\frac{n-1}{2}, 0,\ldots,0\bigg),
\end{equation*}
where there are $k$ entries of $(n-1)/2$ and $k$ entries of $0$ in the Meijer $G$-function.} An exact formula for the PDF of the product of three or more correlated normal random variables is not available in the literature.

\subsection{Cumulative distribution function}\label{sec2.2}

For general parameter values $n\geq1$, $-1<\rho<1$ and $\sigma_X,\sigma_Y>0$, a closed-form formula for the cumulative distribution function (CDF) of $\overline{Z}_n$ is not available. We record some cases for which exact formulas are available. Let $F_{\overline{Z}_n}(x)=\mathbb{P}(\overline{Z}_n\leq x)$. 

Suppose $\rho=0$. Then, by the symmetry of the PDF of $\overline{Z}_n$, it follows that its CDF is given by $F_{\overline{Z}_n}(x)=1/2+\mathrm{sgn}(x)\int_0^{|x|}f_{\overline{Z}_n}(t)\,\mathrm{d}t$, where $\mathrm{sgn}(x)$ denotes the sign of $x$. On calculating the integral using (\ref{besint}), we have that, for $x\in\mathbb{R}$,
\begin{align*}F_{\overline{Z}_n}(x)=\frac{1}{2}+\frac{x}{2s_n}\bigg[K_{\frac{n-1}{2}}\bigg(\frac{|x|}{s_n}\bigg)\mathbf{L}_{\frac{n-3}{2}}\bigg(\frac{|x|}{s_n}\bigg)+\mathbf{L}_{\frac{n-1}{2}}\bigg(\frac{|x|}{s_n}\bigg)K_{\frac{n-3}{2}}\bigg(\frac{|x|}{s_n}\bigg)\bigg],
\end{align*}
where $\mathbf{L}_\nu(x)$ is a modified Struve function of the first kind.
By using the fact that $\overline{Z}_n$ follows the McKay Type II distribution (a VG distribution), other formulas for the CDF in the case $\rho=0$ can be obtained from results of Jankov Ma\v{s}irevi\'c and Pog\'any \cite{jp} and Nadarajah, Srivastava and Gupta \cite{nsg07}.

Now suppose that $n\in 2\mathbb{Z}^+$ and $-1<\rho<1$. Then, making use of the formula (\ref{nsimp}) for the PDF of $\overline{Z}_n$ and calculating $F_{\overline{Z}_n}(x)=\int_{-\infty}^xf_{\overline{Z}_n}(t)\,\mathrm{d}t$ for $x\leq0$, and $F_{\overline{Z}_n}(x)=1-\int_{x}^\infty f_{\overline{Z}_n}(t)\,\mathrm{d}t$ for $x>0$, yields the following formulas. For $x\leq0$,
\begin{align*}F_{\overline{Z}_n}(x)=\frac{(1-\rho)^{n/2}}{2^{n/2}(n/2-1)!}\sum_{j=0}^{n/2-1}\frac{(n/2-1+j)!}{(n/2-1-j)!j!}\bigg(\frac{1+\rho}{2}\bigg)^j\Gamma\bigg(\frac{n}{2}-j,\frac{-x}{s_n(1-\rho)}\bigg),
\end{align*}
and, for $x>0$,
\begin{align*}F_{\overline{Z}_n}(x)=1-\frac{(1+\rho)^{n/2}}{2^{n/2}(n/2-1)!}\sum_{j=0}^{n/2-1}\frac{(n/2-1+j)!}{(n/2-1-j)!j!}\bigg(\frac{1-\rho}{2}\bigg)^j\Gamma\bigg(\frac{n}{2}-j,\frac{x}{s_n(1+\rho)}\bigg),
\end{align*}
where $\Gamma(a,x)=\int_x^\infty t^{a-1}\mathrm{e}^{-t}\,\mathrm{d}t$ is the upper incomplete gamma function. 
Nadarajah, Srivastava and Gupta \cite{nsg07} also gave an analogue of these formulas for the case $\rho=0$ for the CDF of the McKay Type II distribution.

Suppose now that we are in the general setting $n\geq1$, $-1<\rho<1$ and $\sigma_X,\sigma_Y>0$. As a closed-form formula is not available for the CDF, the following asymptotic approximations for the tail probabilities are of interest. Let $\bar{F}_{\overline{Z}_n}(x)=1-F_{\overline{Z}_n}(x)=\mathbb{P}(\overline{Z}_n> x)$. Then, from (\ref{pdf2}) and the limiting form (\ref{kintap}), we get that
\begin{align}\label{barbar}\bar{F}_{\overline{Z}_n}(x)\sim\frac{(1+\rho)}{2^{n/2}s_n^{n/2-1}\Gamma(n/2)}x^{n/2-1}\exp\bigg(-\frac{x}{s_n(1+\rho)}\bigg), \quad x\rightarrow\infty,
\end{align}
and, by symmetry,
\begin{align*}F_{\overline{Z}_n}(x)\sim\frac{(1-\rho)}{2^{n/2}s_n^{n/2-1}\Gamma(n/2)}(-x)^{n/2-1}\exp\bigg(\frac{x}{s_n(1-\rho)}\bigg), \quad x\rightarrow-\infty.
\end{align*}
 
Upper and lower bounds on $F_{\overline{Z}_n}(x)$ and $\overline{F}_{\overline{Z}_n}(x)$ can also be obtained by using bounds for the integral $\int_x^\infty\mathrm{e}^{\beta t}t^\nu K_{\nu}(t)\,\mathrm{d}t$, $x>0$, $-1<\beta<1$, $\nu>-1/2$, given in Gaunt \cite{gaunt ineq1, gaunt ineq3}. \gr{As an example, inequality (2.10) of Gaunt \cite{gaunt ineq1} states that $\int_x^\infty\mathrm{e}^{\beta t}t^\nu K_{\nu}(t)\,\mathrm{d}t< (1-\beta)^{-1}\mathrm{e}^{\beta x}x^\nu K_\nu(x)$, for $x>0$, $0\leq\beta<1$, $-1/2<\nu<1/2$. Applying this bound with $\beta=\rho$ and $\nu=0$ gives that, for $0\leq\rho<1$, $x>0$,
\begin{align}\bar{F}_{Z}(x)&=\int_{x}^\infty \frac{1}{\pi s\sqrt{1-\rho^2}}\exp\bigg(\frac{\rho t}{s(1-\rho^2)}\bigg)K_0\bigg(\frac{t}{s(1-\rho^2)}\bigg)\,\mathrm{d}t\nonumber\\
&=\int_{x/(s(1-\rho^2))}^\infty \frac{\sqrt{1-\rho^2}}{\pi}\mathrm{e}^{\rho u}K_0(u)\,\mathrm{d}u\nonumber\\
\label{cumbd}&< \frac{1}{\pi}\sqrt{\frac{1+\rho}{1-\rho}}\exp\bigg(\frac{\rho x}{s(1-\rho^2)}\bigg)K_0\bigg(\frac{x}{s(1-\rho^2)}\bigg)=s(1+\rho)f_Z(x).
\end{align} 
Applying the limiting form (\ref{Ktendinfinity}) to the upper bound (\ref{cumbd}) and comparing to (\ref{barbar}) shows that the upper bound (\ref{cumbd}) is tight as $x\rightarrow\infty$.
}

\subsection{Generating functions and infinite divisibility}\label{sec2.4}

The moment generating function of $Z$ was obtained by Craig \cite{craig} 
and the moment generating (and characteristic) function of $\overline{Z}_n$ can then be deduced from basic properties of moment generating (and characteristic) functions (see Nadarajah and Pog\'{a}ny \cite{np16}):
\begin{align*}M_{\overline{Z}_n}(t)=\mathbb{E}[\mathrm{e}^{t\overline{Z}_n}]=\big(1-2\rho s_nt-s_n^2(1-\rho^2)t^2\big)^{-n/2},
\end{align*}
which exists if $-1/(1-\rho)<s_n t<1/(1+\rho)$ (see (\ref{tail1}) and (\ref{tail2})). The characteristic function is
\begin{equation}\label{cfcf}\varphi_{\overline{Z}_n}(t)=\mathbb{E}[\mathrm{e}^{\mathrm{i}t\overline{Z}_n}]=\big(1-2\rho s_n\mathrm{i}t+s_n^2(1-\rho^2)t^2\big)^{-n/2}.
\end{equation}
Alternatively, given the formula (\ref{pdf2}) for the PDF of $\overline{Z}_n$, the moment generating function is obtained from a simple integration of $M_{\overline{Z}_n}(t)=\int_{-\infty}^\infty\mathrm{e}^{tx}f_{\overline{Z}_n}(x)\,\mathrm{d}x$; see, for example, Eberlein and Hammerstein \cite{eberlein} for a calculation for the more general generalized hyperbolic distribution (which includes the VG distribution as a limiting case). 

The \gr{cumulant} generating function $K_{\overline{Z}_n}(t)=\log\mathbb{E}[\mathrm{e}^{t\overline{Z}_n}]$, defined for $-1/(1-\rho)<s_n t<1/(1+\rho)$, is given by
\begin{align}\label{kfkf}K_{\overline{Z}_n}(t)&=-\frac{n}{2}\log\big(1-2\rho s_nt-s_n^2(1-\rho^2)t^2\big) \nonumber \\
&=-\frac{n}{2}\log\big(1-s_n(\rho+1)t\big)-\frac{n}{2}\log\big(1-s_n(\rho-1)t\big).
\end{align}

\gr{The distribution of $\overline{Z}_n$ is infinitely divisible. This is easily seen, because  $\overline{Z}_n$ follows a VG distribution, and the VG distribution is a special case of the generalized hyperbolic distribution which is infinitely divisible (see Barndorff-Nielsen and Halgreen \cite{bh77}). That the distribution of $\overline{Z}_n$ is infinitely divisible can also be inferred from the representation of $\overline{Z}_n$ given in part 3 of Section \ref{sec2.5}}.

\subsection{Representation in terms of other random variables}\label{sec2.5}

\noindent{1.} Let $X_1,\ldots,X_n$ and $W_1,\ldots,W_n$ be independent $N(0,1)$ random variables. Then
\begin{equation*}\overline{Z}_n=_d s_n\sum_{j=1}^n\big(\sqrt{1-\rho^2}X_jW_j+\rho X_j^2\big).
\end{equation*}
It suffices to show that $Z=_d s(\sqrt{1-\rho^2}X_1W_1+\rho X_1^2\gr{)}$\gr{; this was done in Gaunt \cite{gaunt prod} and we repeat the simple steps here.} For ease of notation, we suppose that $s=1$, with the general case following by rescaling. It is straightforward to verify that $X$ and $W=(Y-\rho X)/\sqrt{1-\rho^2}$ are independent $N(0,1)$ random variables. Thus, $Z=XY=X(\sqrt{1-\rho^2}W+\rho X)=\sqrt{1-\rho^2}XW+\rho X^2$, as required. 

\vspace{3mm}

\noindent{2.} Suppose that $S\sim \chi_{(n)}^2$ and $T\sim N(0,1)$ are independent. Then
\begin{equation}\label{rep2}\overline{Z}_n=_d \rho s_n S+s_n\sqrt{1-\rho^2}\sqrt{S}T.
\end{equation}
The representation (\ref{rep2}) follows from the representation (\ref{vgrep}) of $\overline{Z}_n$ as a VG random variable and Proposition 4.1.2 of Kotz, Kozubowski and Podg\'{o}rski \cite{kkp01}, which states that if $S'\sim \Gamma(r/2,1/2)$ and $T'\sim N(0,1)$ are independent, then $\theta S'+\sigma \sqrt{S'}T'\sim \mathrm{VG}(r,\theta,\sigma,0)$.

\vspace{3mm}

\gr{\noindent{3.} For any $m\geq1$, we can write $\overline{Z}_n=_dS_1+S_2+\cdots+S_m$, where $S_1,S_2,\ldots,S_m$ are independent $\mathrm{VG}(n/m,\rho s_n,s_n\sqrt{1-\rho^2},0)$ random variables. 
This follows from the representation (\ref{vgrep}) of $\overline{Z}_n$ as a VG distribution and the convolution property that if $V_1\sim\mathrm{VG}(r_1,\theta,\sigma,0)$ and $V_2\sim\mathrm{VG}(r_2,\theta,\sigma,0)$ are independent then $V_1+V_2\sim\mathrm{VG}(r_1+r_2,\theta,\sigma,0)$ (Bibby and S{\o}rensen \cite{bibby}), which is easily verified by using the fact that $\varphi_V(t)=(1-2\theta t+\sigma^2t^2)^{-r/2}$ is the characteristic function of $V\sim\mathrm{VG}(r,\theta,\sigma,0)$.}

\vspace{3mm}

\noindent{4.} Suppose that $V$ and $V'$ are independent $\chi_{(n)}^2$ random variables. Then
\begin{equation}\label{gamrep}\overline{Z}_n=_d\frac{s_n}{2}(1+\rho)V-\frac{s_n}{2}(1-\rho)V'.
\end{equation}
This representation of $\overline{Z}_n$ is easily proved using a standard characteristic function argument with the formula (\ref{cfcf}) for the characteristic function of $\overline{Z}_n$ and the standard formula for the characteristic function of the $\chi_{(n)}^2$ distribution, $\varphi_V(t)=(1-2\mathrm{i}t)^{-n/2}$.

From the representation (\ref{gamrep}) and the representation of the $\chi_{(n)}^2$ distribution as the sum of the squares of $n$ independent $N(0,1)$ random variables, it follows that 
\[\overline{Z}_n=_d\rho s+\sum_{j=1}^{2n}\lambda_j(N_j^2-1),\]
where $N_1,\ldots,N_{2n}$ are independent $N(0,1)$ random variables, and $\lambda_1=\ldots=\lambda_n=s_n(1+\rho)/2$ and $\lambda_{n+1}=\ldots=\lambda_{2n}=s_n(\rho-1)/2$. Thus, $\overline{Z}_n$ is a member of the second Wiener chaos (see Section 2.2 of Nourdin and Peccati \cite{np12}).

\vspace{3mm}

\noindent{5.} If $n\in 2\mathbb{Z}^+$, then the independent random variables $V\sim\chi_{(n)}^2$ and $V'\sim\chi_{(n)}^2$ can be expressed in terms of independent uniform $U(0,1)$ random variables $U_1,\ldots,U_{n}$ as $V=_d\sum_{j=1}^{n/2}\log U_j$ and $V'=_d\sum_{j=n/2+1}^{n}\log U_j$. Therefore, from (\ref{gamrep}), 
\begin{equation}\label{unifrep}\overline{Z}_n=_d\frac{s_n}{2}(1+\rho)\sum_{j=1}^{n/2}\log U_j-\frac{s_n}{2}(1-\rho)\sum_{j=n/2+1}^{n}\log U_j.
\end{equation}
The representation (\ref{unifrep}) is convenient for simulating the distribution of $\overline{Z}_n$ when $n\in 2\mathbb{Z}^+$. If $n$ is an odd integer, the distribution of  $\overline{Z}_n$ can be simulated using the representation (\ref{gamrep}) and simulating the chi-square distributions of $V$ and $V'$ using methods for simulating gamma distributions, as given in Chapter 9, Section 3 of Devroye \cite{d86}.

\subsection{Stein characterisation}\label{sec2.6}

Let $W$ be a real-valued random variable. Then $W=_d\overline{Z}_n$ if and only if
\begin{align}\label{char1}\mathbb{E}\big[s_n^2(1-\rho^2)Wg''(W)+(ns_n^2(1-\rho^2)+2\rho s_n W)g'(W)
+(\rho s-W)g(W)\big]&=0
\end{align}
for all twice differentiable $g:\mathbb{R}\rightarrow\mathbb{R}$ such that the expectations $\mathbb{E}|g(\overline{Z})|$, $\mathbb{E}|\overline{Z}g(\overline{Z})|$, $\mathbb{E}|g'(\overline{Z})|$, $\mathbb{E}|\overline{Z}g'(\overline{Z})|$ and $\mathbb{E}|\overline{Z}g''(\overline{Z})|$ are all finite. 

Necessity of the Stein characterisation (\ref{char1}) was established by Gaunt \cite{gaunt prod2} (here we have corrected a typo in the characterising equation (\ref{char1})). Let us give a simple proof of sufficiency. Suppose $W$ is a real-valued random variable. To simplify notation, we set $s=1$, with the general case following by rescaling. Taking $g(x)=\mathrm{e}^{\mathrm{i}tx}$ (which is twice differentiable and bounded) in (\ref{char1}) and setting $\varphi(t)=\mathbb{E}[\mathrm{e}^{\mathrm{i}tW}]$ yields the differential equation
\begin{equation}\label{odee}\bigg(\frac{(1-\rho^2)}{n^2}t^2-\frac{2\rho\mathrm{i}}{n}t+1\bigg)\varphi'(t)+\bigg(\frac{(1-\rho^2)}{n}t-\rho\mathrm{i}\bigg)\varphi(t)=0.
\end{equation}
Note that $g(x) = \mathrm{e}^{\mathrm{i}tx}$ is a complex-valued function; here we applied the characterising equation to the real and imaginary parts of $g$, which are real-valued functions. Solving (\ref{odee}) subject to the condition $\varphi(0)=1$ gives that $\varphi(t)=(1-2\rho\mathrm{i}t/n+(1-\rho^2)t^2/n^2)^{-n/2}$, which is the characteristic function (\ref{cfcf}) of $\overline{Z}_n$. Thus, $W=_d\overline{Z}_n$, which proves sufficiency.

Stein characterisations are most commonly used as part of Stein's method (Stein \cite{stein}) to prove quantitative limit theorems in probability theory; however, they can also be used to derive distributional properties. For example, setting $g_1(x)=x^k$ and $g_2(x)=(x-\mathbb{E}[\overline{Z}_n])^k=(x-\rho s)^k$ in (\ref{char1}) yields the following recursions for the $k$-th raw moment $\mu_k'=\mathbb{E}[\overline{Z}_n^k]$ and the $k$-th central moment $\mu_k=\mathbb{E}[(\overline{Z}_n-\mathbb{E}[\overline{Z}_n])^k]$:
\begin{align}\label{rec1}\mu_{k+1}'=(n+2k)\rho s_n\mu_k'+k(n+k-1)s_n^2(1-\rho^2)\mu_{k-1}', \quad k\geq1,
\end{align}
and
\begin{align}\mu_{k+1}&=2k\rho s_n\mu_k+ks_n^2\big(n+k-1+(n-k+1)\rho^2\big)\mu_{k-1}\nonumber\\
\label{rec2}&\quad+k(k-1)ns_n^3\rho(1-\rho^2)\mu_{k-2}, \quad k\geq2.
\end{align}
These recurrences allow for the efficient computation of lower order raw (central) moments given just the first raw (first and second central) moments.

\subsection{Moments and cumulants}\label{sec2.7}

The mean and variance of $\gr{\overline{Z}_n}$ can be easily calculated through several approaches; for example, using the standard method of extracting lower order moments via moment generating functions.
We have
\begin{align}\label{mean}\mathbb{E}[\gr{\overline{Z}_n}]&=\rho s, \\
\label{var}\mathrm{Var}(\gr{\overline{Z}_n})&=ns_n^2(1+\rho^2).
\end{align}
Lower order raw and central moments are then readily obtained from the recurrences (\ref{rec1}) and (\ref{rec2}). The first four raw moments are given by
\begin{align*}\mu_1'&=\rho s, \\
\mu_2'&=ns_n^2\big(1+(n+1)\rho^2\big), \\
\mu_3'&=n\rho s_n^3\big(3(n+2)+(n+1)(n+2)\rho^2\big), \\
\mu_4'&=ns_n^4\big(3(n+2)+6(n+2)(n+3)\rho^2+(n+1)(n+2)(n+3)\rho^4\big),
\end{align*}
and the first four central moments are given by
\begin{align*}\mu_1&=0,\\
\mu_2&=ns_n^2(1+\rho^2),  \\
\mu_3&=2n\rho s_n^3(3+\rho^2), \\
\mu_4&=3ns_n^4\big((n+2)+2(n+6)\rho^2+(n+2)\rho^4\big).
\end{align*}
We thus deduce that the skewness $\gamma_1=\mu_3/\mu_2^{3/2}$ and kurtosis $\beta_2=\mu_4/\mu_2^2$ are given by
\begin{align*}\gamma_1=\frac{2\rho(3+\rho^2)}{\sqrt{n}(1+\rho^2)^{3/2}}, \quad
\beta_2=\frac{3(n+2)+6(n+6)\rho^2+3(n+2)\rho^4}{n(1+\rho^2)^2},
\end{align*}
whilst the excess kurtosis $\gamma_2=\mu_4/\mu_2^2-3$ is
\[\gamma_2=\frac{6+36\rho^2+6\rho^4}{n(1+\rho^2)^2}.\]
\gr{Formulas} for the skewness \gr{and kurtosis} of the product of two correlated normal variables with possibly non-zero means \gr{are} given by \gr{Seijas-Mac\'{i}as et al.\ \cite{sym2020} and} Ware and Lad \cite{ware}. \gr{Formulas for the first four central moments of the product of two correlated normal variables with possibly non-zero means are given by Haldane \cite{hh1}.}

\begin{figure}[h!]
\centering
\includegraphics[width=6in]{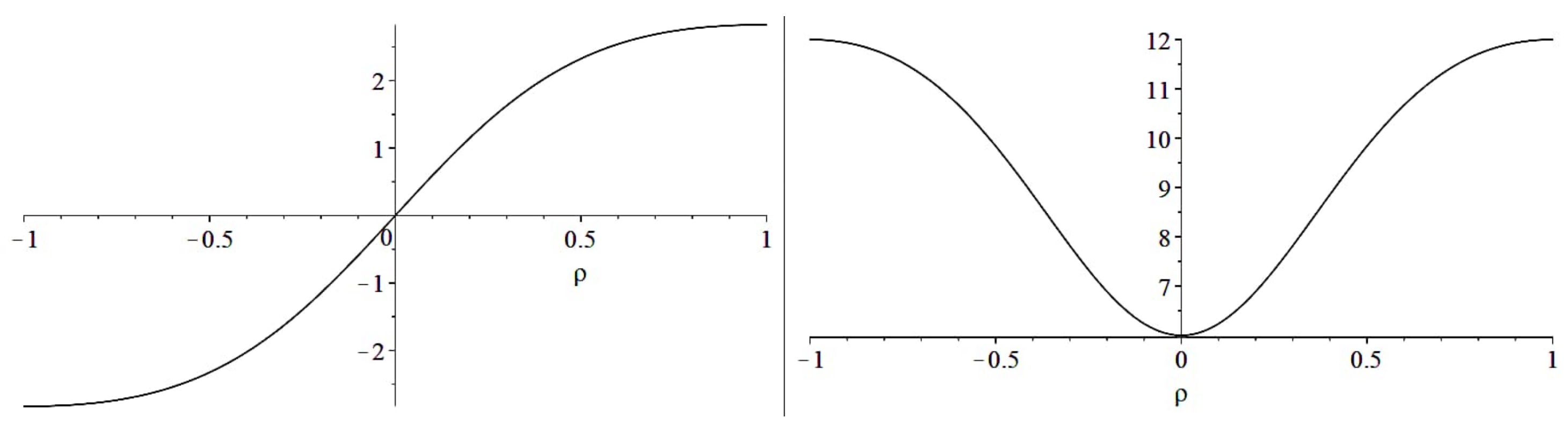}
\caption{Skewness $\gamma_1$ (left panel) and excess kurtosis $\gamma_2$ (right panel) with $n=1$}
\end{figure}

Higher order moments are easily calculated using the representation (\ref{gamrep}) of $\overline{Z}_n$ and the moment formula $\mathbb{E}[V^k]=2^k\Gamma(n/2+k)/\Gamma(n/2)=n(n+2)(n+4)\cdots(n+2k-2)$, $k\geq1$, for $V\sim\chi_{(n)}^2$. For $k\geq1$,
\begin{align*}\mu_k'&=\frac{s_n^k}{2^k}\sum_{j=0}^k\binom{k}{j}(\rho+1)^j(\rho-1)^{k-j}\mathbb{E}[V^j]\mathbb{E}[V^{k-j}] \\
&=\frac{s_n^k}{(\Gamma(n/2))^2}\sum_{j=0}^k\binom{k}{j}\Gamma\Big(\frac{n}{2}+j\Big)\Gamma\Big(\frac{n}{2}+k-j\Big)(\rho+1)^j(\rho-1)^{k-j}.
\end{align*}
In the case $n=1$, the following formula is available (see Kan \cite{k08}):
\begin{equation*}\mu_k'=\begin{cases} \displaystyle \frac{s^k(k!)^2}{2^k}\sum_{j=0}^{k/2}\frac{(2\rho)^{2j}}{((k/2-j)!)^2(2j)!}, & \quad \text{if $k$ is even}, \\
\displaystyle \frac{s^k(k!)^2}{2^k}\sum_{j=0}^{(k-1)/2}\frac{(2\rho)^{2j+1}}{((k/2-1/2-j)!)^2(2j+1)!}, & \quad \text{if $k$ is odd}. \end{cases}
\end{equation*}
Returning to the general case $n\geq1$, the raw moments of $\overline{Z}_n$ can be expressed in terms of the hypergeometric function (see Chapter 15 of Olver et al.\ \cite{olver}).
For $k\geq1$,
\begin{align*}\mu_k'&=\frac{s_n^k(n+k-1)!k!}{(1-\rho^2)^{n/2}\Gamma(n/2+k+1)\Gamma(n/2)}\times\\
&\quad\times\bigg[(-1)^k(1-\rho)^{n+k}\,{}_2F_1\bigg(n+k,\frac{n}{2};\frac{n}{2}+k+1;-\frac{1-\rho}{1+\rho}\bigg)\\
&\quad\quad\quad+(1+\rho)^{n+k}\,{}_2F_1\bigg(n+k,\frac{n}{2};\frac{n}{2}+k+1;-\frac{1+\rho}{1-\rho}\bigg)\bigg].
\end{align*}
This formula is a re-parameterisation (using (\ref{vgrep})) of the result of Theorem 2 of Holm and Alouini \cite{ha04}, which gives formulas for the raw moments of the Mckay Type II distribution (a VG distribution). The proof of Holm and Alouini \cite{ha04} involves writing $\mu_k'=\int_{-\infty}^\infty x^kf_{\overline{Z}_n}(x)\,\mathrm{d}x=\int_{-\infty}^0 x^kf_{\overline{Z}_n}(x)\,\mathrm{d}x+\int_{0}^\infty x^kf_{\overline{Z}_n}(x)\,\mathrm{d}x=(-1)^k\int_{0}^\infty x^kf_{\overline{Z}_n}(-x)\,\mathrm{d}x+\int_{0}^\infty x^kf_{\overline{Z}_n}(x)\,\mathrm{d}x$, calculating the latter two integrals using (\ref{hyint}) and then simplifying. The absolute moments can be calculated similarly by writing $\mathbb{E}|\overline{Z}_n^k|=\int_{0}^\infty x^kf_{\overline{Z}_n}(-x)\,\mathrm{d}x+\int_{0}^\infty x^kf_{\overline{Z}_n}(x)\,\mathrm{d}x$. This gives that, for $k\geq1$,
\begin{align*}\mathbb{E}|\overline{Z}_n^k|&=\frac{s_n^kn+k-1)!k!}{(1-\rho^2)^{n/2}\Gamma(n/2+k+1)\Gamma(n/2)}\times\\
&\quad\times\bigg[(1-\rho)^{n+k}\,{}_2F_1\bigg(n+k,\frac{n}{2};\frac{n}{2}+k+1;-\frac{1-\rho}{1+\rho}\bigg)\\
&\quad\quad\quad+(1+\rho)^{n+k}\,{}_2F_1\bigg(n+k,\frac{n}{2};\frac{n}{2}+k+1;-\frac{1+\rho}{1-\rho}\bigg)\bigg].
\end{align*}

Higher order central moments of $\overline{Z}_n$ can be calculated using the representation (\ref{gamrep}) of $\overline{Z}_n$ and the formula $\mathbb{E}[(V-\mathbb{E}[V])^k]=2^kU(-k,1-k-n/2,-n/2)$, where $U(a,b,x)$ is a confluent hypergeometric function of the second kind (see Weisstein \cite{mathworld}). We note that in the case $a=-m$, $m=0,1,2,\ldots$, the function $U(a,b,x)$ is a polynomial:
$U(-m,b,x)=(-1)^m\sum_{j=0}^m\binom{m}{j}(b+j)_{m-j}(-x)^j$, and the Pochhammer symbol is $(a)_0=0$ and $(a)_j=a(a+1)(a+2)\cdots(a+j-1)$, $j\geq1$ (see Section 13.2(i) of Olver et al.\ \cite{olver}). We have, for $k\geq1$,
\begin{align*}\mu_k&=\frac{s_n^k}{2^k}\sum_{j=0}^k\binom{k}{j}(\rho+1)^j(\rho-1)^{k-j}\mathbb{E}[(V-\mathbb{E}[V])^j]\mathbb{E}[(V-\mathbb{E}[V])^{k-j}] \\
&=s_n^k\sum_{j=0}^k\binom{k}{j}U\bigg(-j,1-j-\frac{n}{2},-\frac{n}{2}\bigg)U\bigg(j-k,1+j-k-\frac{n}{2},-\frac{n}{2}\bigg)\times\\
&\quad\times(1+\rho)^j(\rho-1)^{k-j}.
\end{align*}   

Suppose now that $\rho=0$. Then, from (\ref{rep2}), $\overline{Z}_n=_d s_n\sqrt{S}T$, where $S\sim\chi_{(n)}^2$ and $T\sim N(0,1)$ are independent. If $k\geq1$ is odd then $\mathbb{E}[T^k]=0$, whilst for any $\alpha>0$, $\mathbb{E}|T^\alpha|=\pi^{-1/2}2^{\alpha/2}\Gamma((\alpha+1)/2)$ and $\mathbb{E}|S^\alpha|=\Gamma(n/2+\alpha)/\Gamma(n/2)$. Therefore, for $k\geq1$,
\begin{equation*}\mu_k'=\begin{cases} \displaystyle \frac{2^{k/2}s_n^k}{\sqrt{\pi}}\frac{\Gamma((n+k)/2)\Gamma((k+1)/2)}{\Gamma(n/2)}, & \quad \text{if $k$ is even}, \\
\displaystyle 0, & \quad \text{if $k$ is odd}. \end{cases}
\end{equation*}

Suppose now that $-1<\rho<1$ and $n\geq1$. The cumulants of $\overline{Z}_n$ are readily calculated by either Taylor expanding the logarithms in (\ref{kfkf}), or by working with the representation (\ref{gamrep}) of $\overline{Z}_n$. In the latter approach, one uses the standard properties that for independent random variables $S_1$ and $S_2$ and constant $c$, the $k$-th cumulant satisfies $\kappa_k(cS_1)=c^k\kappa_k(S_1)$ and $\kappa_k(S_1+S_2)=\kappa_k(S_1)+\kappa_k(S_2)$, and $\kappa_k(V)=2^{k-1}(k-1)!\,n$ for $V\sim\chi_{(n)}^2$. For $k\geq1$, 
\begin{align*}\kappa_k=\frac{n s_n^k}{2}(k-1)!\big[(1+\rho)^k+(\rho-1)^k\big].
\end{align*}
 In particular,
\begin{align*}\kappa_1&=\rho s, \\
\kappa_2&=ns_n^2(1+\rho^2), \\
\kappa_3&=2\rho n s_n^3(3+\rho^2), \\
\kappa_4&=6ns_n^4(1+6\rho^2+\rho^4),\\
\gr{\kappa_5}&\gr{=24\rho ns_n^5(5+10\rho^2+\rho^4),} \\
\gr{\kappa_6}&\gr{=120n s_n^6(1+\rho^2)(1+14\rho^2+\rho^4)}.
\end{align*}
\gr{A general formula for the cumulants of the product of two correlated normal variables with possibly non-zero means is given by Craig \cite{craig}.}

\subsection{Mode and median}\label{sec2.8}

A detailed study of the mode and median of the generalized hyperbolic and VG distributions was recently carried out by Gaunt and Merkle \cite{gm21}; here we provide a synthesis of this work in the case of the product of correlated zero mean normal random variables. 

The distribution of $\overline{Z}_n$ is unimodal. This follows because $\overline{Z}_n$ is self-decomposable (it is infinitely divisible) and self-decomposable distributions are unimodal (Yamazato \cite{y78}). Let $M_n$ denote the mode of $\overline{Z}_n$. It is clear from (\ref{pmutend}) and (\ref{lapd}) that 
\begin{equation}\label{mode0}M_1=M_2=0.
\end{equation}
 Suppose now that $n\geq3$. Then, by applying the differentiation formula (\ref{ddbk}) to the PDF (\ref{pdf2}), we deduce that $M_n=\mathrm{sgn}(\rho)\cdot x^*$, where $x^*$ is the unique positive solution of the equation 
\begin{equation}\label{xstar}K_{\frac{n-3}{2}}\bigg(\frac{x}{s_n(1-\rho^2)}\bigg)=|\rho|K_{\frac{n-1}{2}}\bigg(\frac{x}{s_n(1-\rho^2)}\bigg).
\end{equation}
In the cases $n=4$ and $n=6$, we can apply the formulas in (\ref{special2}) to (\ref{xstar}) to obtain simple algebraic equations for $x^*$ which when solved yield the exact expressions
\begin{equation*}M_4=\frac{\rho s}{4}(1+|\rho|), \quad M_6 =\frac{\rho s}{12}(1+|\rho|)\bigg(3-\frac{1}{|\rho|}+\sqrt{\frac{1}{\rho^2}+\frac{6}{|\rho|}-3}\bigg).
\end{equation*}
For $n=8$ and $n=10$, we can apply the formula (\ref{special}) to (\ref{xstar}) to obtain cubic and quartic equations for $x^*$, although the solutions to these equations are too complicated to be worth reporting. For other values of $n$, exact formulas for $M_n$ are not available, except for the case $\rho=0$ in which case the mode is equal to zero.

\begin{figure}[h!]
\centering
\includegraphics[width=5.0in]{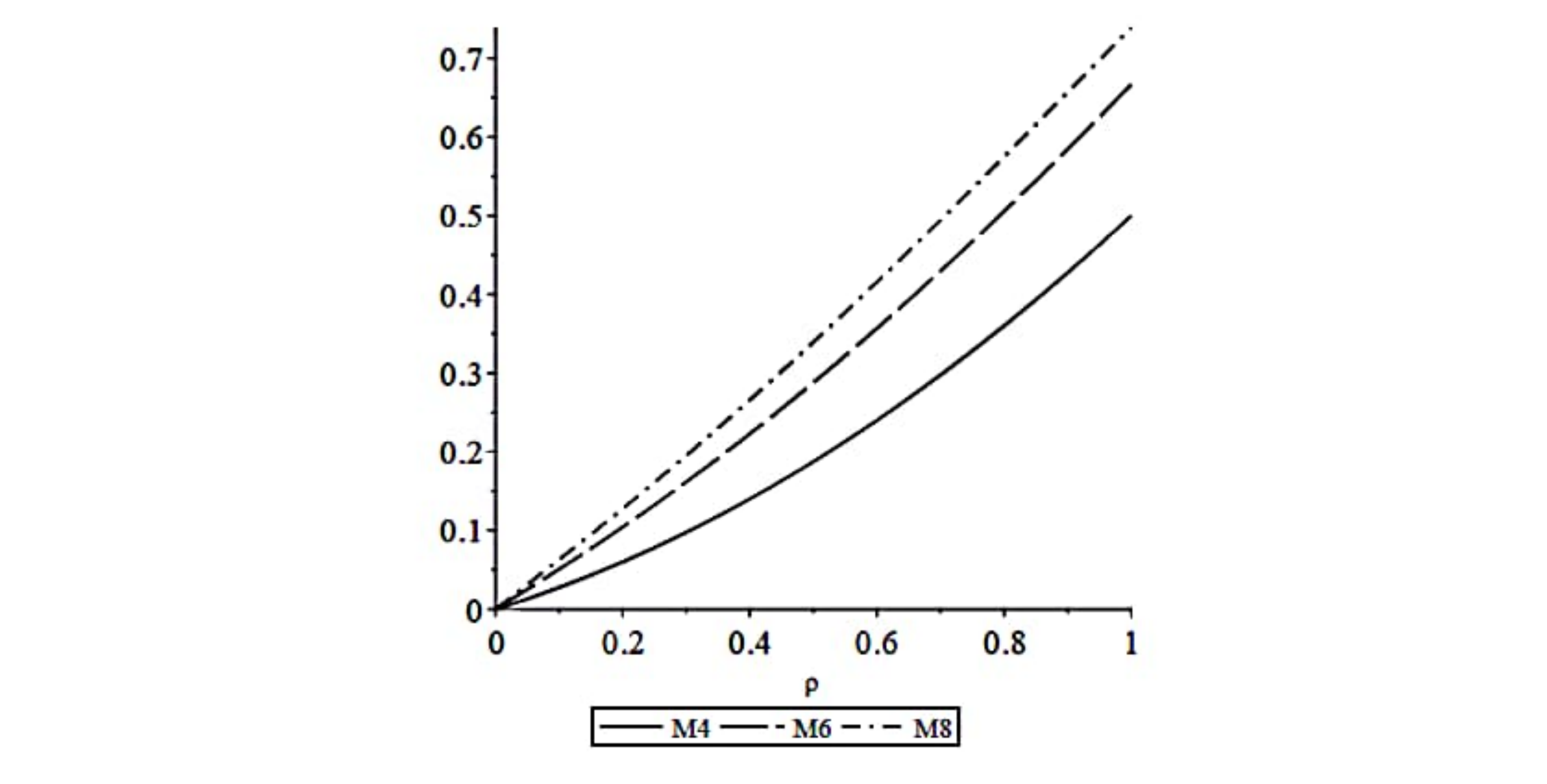}
\caption{Plots of $M_4$, $M_6$ and $M_8$ when $s=1$}
\end{figure}

Whilst exact formulas are not available for $M_n$ for general $n\geq3$, simple and accurate upper and lower bounds can be derived.
In what follows, we fix $\rho>0$; bounds for the case $\rho<0$ follow by symmetry. Applying the lower and upper bounds of (\ref{seg1}) to (\ref{xstar}) leads to simple algebraic equations for $x^*$ which when solved lead to the upper and lower bounds of the following two-sided inequality, respectively:
\begin{equation}\label{mode1}\rho s(1-3/n)<M_n<\rho s(1-2/n), \quad n\geq3.
\end{equation}
There is equality in the upper bound when $n=2$ (see (\ref{mode0})). Applying inequality (\ref{seg3}) to (\ref{xstar}) yields that, for $n\geq4$, 
\begin{equation}\label{mode2}M_n\geq\frac{\rho s}{2}\bigg[1-\frac{2}{n}+\sqrt{\rho^2\Big(1-\frac{2}{n}\Big)^2+(1-\rho^2)\Big(1-\frac{4}{n}\Big)^2}\bigg],
\end{equation}
with equality if and only if $n=4$. Inequality (\ref{mode2}) is more accurate than the lower bound in (\ref{mode1}) for $n\geq4$. The bounds in (\ref{mode1}) and (\ref{mode2}) were derived by Gaunt and Merkle \cite{gm21} by similar considerations.


If $\rho=0$, the median of $\overline{Z}_n$ is equal to zero. The median of $\overline{Z}_2$ is given by
\begin{equation}\label{medform}\mathrm{Med}(\overline{Z}_2)=\begin{cases} \displaystyle\frac{s}{2}(1+\rho)\log(1+\rho), & \quad \rho\geq0, \\
\displaystyle-\frac{s}{2}(1-\rho)\log(1-\rho), & \quad \rho<0, \end{cases}
\end{equation}
which follows from the well-known formula for the  median of the asymmetric Laplace distribution (see Kozubowski and Podg\'{o}rski \cite{kp01}). Otherwise, an exact closed-form formula is not available for $\mathrm{Med}(\overline{Z}_n)$, $n\geq1$. Moreover, unlike for the mode, accurate upper and lower bounds for the median have yet to be worked out in the literature. Gaunt and Merkle \cite{gm21} have, however, conjectured accurate bounds for the median of the VG distribution. As $\overline{Z}_n$ follows the VG distribution (see (\ref{vgrep})), we can present conjectured bounds for $\mathrm{Med}(\overline{Z}_n)$. The numerical results of Table \ref{table1} support the conjectured bounds\gr{, in that each entry in the table lies between the conjectured lower and upper bounds}. The results were obtained with \emph{Mathematica} via the same simple numerical procedure for numerically computating medians as described in pp.\ 16--17 of Gaunt and Merkle \cite{gm21}.

\begin{conjecture}\label{conj3}Suppose $\rho>0$. Then it is conjectured that
\begin{equation*}\bigg(1-\frac{1}{n}\bigg)\rho s<\mathrm{Med}(\overline{Z}_n)<\rho s \mathrm{e}^{-2/3n}<\bigg(1-\frac{2}{3n}+\frac{2}{9n^2}\bigg)\rho s, \quad n\geq1,
\end{equation*}
and
\begin{equation}\label{concon}\mathrm{Med}(\overline{Z}_n)\leq\bigg(1-\frac{2(1-\log2)}{n}\bigg)\rho s, \quad n\geq2.
\end{equation}
\end{conjecture}

\begin{table}[h]
\centering
\caption{\normalsize{Median of $\overline{Z}_n$ with $s=1$.}}
\label{table1}
{\normalsize
\begin{tabular}{|c|rrrrr|}
\hline
 \backslashbox{$n$}{$\rho$}       &    0.1 &    0.3 &  0.5 & 0.7 &    0.9    \\
 \hline
1  & 0.0198 & 0.0813 & 0.164 & 0.265 & 0.386  \\
3  & 0.0674 & 0.210 & 0.364 & 0.528  & 0.700  \\
5  & 0.0802 & 0.245 & 0.416 & 0.594  &  0.777   \\
7  & 0.0859 & 0.260 & 0.439 & 0.623 & 0.812    \\ 
10 & 0.0901 & 0.272 & 0.457 & 0.646 & 0.838    \\ 
  \hline
\end{tabular}}
\end{table}

We end by noting that because $\mathbb{E}[\overline{Z}_n]=\rho s$, it follows from (\ref{mode0}) and (\ref{mode1}) that, if $\rho>0$, then $M_n<\mathbb{E}[\overline{Z}_n]$ for all $n\geq1$. In fact, we have the two-sided inequality $2\rho s_n<\mathbb{E}[\overline{Z}_n]-M_n<3\rho s_n$, $n\geq3$, with equality in the lower bound if $n=2$. If the conjectured median bounds hold, then it would follow that, for $\rho>0$, $M_n<\mathrm{Med}(\overline{Z}_n)<\mathbb{E}[\overline{Z}_n]$, meaning that $\overline{Z}_n$ would satisfy the mean-median-mode inequality (Groeneveld and Meeden \cite{gm77} and van Zwet \cite{v79}). We also note that as $n\rightarrow\infty$, $M_n$ and $\mathrm{Med}(\overline{Z}_n)$ converge to $\mathbb{E}[\overline{Z}_n]=\rho s$, which is to be expected given that $\overline{Z}_n$ is approximately normally distributed (for which the mean, median and mode are equal) for large $n$.

\gr{\section{Application as a limiting distribution}\label{sec3}

In Section \ref{sec2.5}, we saw that the distributions of $Z$ and $\overline{Z}_n$ have simple representations in terms of independent standard normal, chi-square and variance-gamma random variables. This feature means that $Z$ and $\overline{Z}_n$ are natural candidates for limiting distributions. 
In this section, we shall see that $Z$ and $\overline{Z}_n$ are limiting distributions for sequences of double Wiener-It\^{o} integrals. Variance-gamma approximations for double Wiener-It\^{o} integrals have been studied by Azmoodeh, Eichelsbacher and Th\"ale \cite{aet21},  Eichelsbacher and Th\"ale \cite{eichelsbacher} and Gaunt \cite{gaunt vgiii}. Here we present a synthesis of the results in the case of the product of correlated mean zero normal random variables.

\subsection{A six moment theorem for double Wiener-It\^{o} integrals}

We first introduce some notation and terminology. Let $\mathfrak{H}$ be a real separable Hilbert space and let $\mathfrak{H}^{\odot 2}$ denote the second symmetric tensor product of $\mathfrak{H}$. For $f\in \mathfrak{H}^{\odot 2}$, the double Wiener-It\^{o} integral is denoted by $I_2(f)$ (see Definition 2.7.1 of Nourdin and Peccati \cite{np12}). Some of the most important properties of multiple  Wiener-It\^{o} integrals can be found in Section 2.7 of Nourdin and Peccati \cite{np12}. If $f\in L^2([0,T]^2,\mathrm{d}t)$ is symmetric then
\[I_2(f)=\int_{[0,T]^2}f(t_1,t_2)\,\mathrm{d}B_{t_1}\,\mathrm{d}B_{t_2},\]
where $B=(B_t)_{t\in[0,T]}$ is a standard two-dimensional Brownian motion (see Exercise 2.7.6 of Nourdin and Peccati \cite{np12}).
Consider also the Wasserstein distance $d_{\mathrm{W}}(F,G)$ and smooth Wasserstein distance $d_{2}(F,G)$ between the distributions of two random elements $F$ and $G$, defined by
\begin{align*}d_{\mathrm{W}}(F,G):&=\sup_{h\in\mathcal{H}_{\mathrm{W}}}|\mathbb{E}h(F)-\mathbb{E}h(G)|, \\
d_{2}(F,G):&=\sup_{h\in\mathcal{H}_2}|\mathbb{E}h(F)-\mathbb{E}h(G)|,
\end{align*}
where
\begin{align*}\mathcal{H}_{\mathrm{W}}&=\{h:\mathbb{R}\rightarrow\mathbb{R}\,|\,\text{$h'$ is Lipschitz, $\|h'\|_\infty\leq1$}\},\\\mathcal{H}_2&=\{h:\mathbb{R}\rightarrow\mathbb{R}\,|\,\text{$h'$ is Lipschitz, $\|h'\|_\infty\leq1$, $\|h''\|_\infty\leq1$}\}. 
\end{align*}
Note that $d_{2}(F,G)\leq d_{\mathrm{W}}(F,G)$ for any random elements $F$ and $G$ for which $d_{\mathrm{W}}(F,G)$ is well-defined. 

Let $F_m=I_2(f_m)$ with $f_m\in \mathfrak{H}^{\odot 2}$, $m\geq1$. Let $\overline{Z}_{n,c}$ be a random variable equal in distribution to $\overline{Z}_{n}-\rho s$, so that $\overline{Z}_{n,c}$ has mean zero. Then recasting Theorem 5.8 of Eichelsbacher and Th\"ale \cite{eichelsbacher} in terms of the distribution of $\overline{Z}_n$, we have that, as $m\rightarrow\infty$, the sequence $(F_m)_{m\geq1}$ converges in distribution to $\overline{Z}_{n,c}$ if and only if $\kappa_i(F_m)\rightarrow\kappa_i(\overline{Z}_{n,c})$, $i=2,3,4,5,6$. The cumulants $\kappa_i(\overline{Z}_{n,c})=\kappa_i(\overline{Z}_{n}-\rho s)$ can be calculated using the formulas of Section \ref{sec2.7} and the standard formula $\kappa_j(S+c)=\kappa_j(S)$, $j\geq2$, for any $c\in\mathbb{R}$. This ``six moment" theorem tells us that the convergence of a sequence of double Wiener-It\^{o} integrals to the distribution of the centered random variable $\overline{Z}_{n,c}$ is determined only by the behaviour of the first six cumulants (equivalently first six moments). This is a product of correlated normal random variables analogue of the celebrated ``fourth moment" theorem for normal approximation of multiple Wiener-It\^{o} integrals of Nualart and Peccati \cite{np05}.

Moreover, quantitative ``six moment" theorems are available. Define 
\[\mathbf{M}(F_m)=\max\{|\kappa_i(F_m)-\kappa_i(\overline{Z}_{n,c})|\,:\,i=2,3,4,5,6\}.\]
Then there exists a constant $C$ only depending on $n$, $\rho$ and $s$ such that
\begin{equation}\label{optvg0}d_{\mathrm{W}}(F_m,\overline{Z}_{n,c})\leq C\sqrt{\mathbf{M}(F_m)}
\end{equation}
(see Eichelsbacher and Th\"ale \cite{eichelsbacher} and Gaunt \cite{gaunt vgiii}), and there exist constants $C_1,C_2>0$ only depending on $n$, $\rho$ and $s$ such that
\begin{equation}\label{optvg}C_1\mathbf{M}(F_m)\leq d_{2}(F_m,\overline{Z}_{n,c})\leq C_2\mathbf{M}(F_m)
\end{equation}
(see Azmoodeh, Eichelsbacher and Th\"ale \cite{aet21}). The upper bound in (\ref{optvg}) improves the bound (\ref{optvg0}) by removing the square root factor, with the improvement coming at the expense of being given with respect to the weaker $d_2$ metric. The rate of convergence in (\ref{optvg}) is optimal, and represents an analogue of the optimal fourth moment theorem for normal approximation of Nourdin and Peccati \cite{np15}.

\subsection{The generalized Rosenblatt process at extreme critical
exponent}

In this section, we see how the bounds (\ref{optvg0}) and (\ref{optvg}) can be applied to obtain bounds on the rate of convergence for one of the main results of Bai and Taqqu \cite{bt17}. Consider the generalized Rosenblatt process $Z_{\gamma_1,\gamma_2}(t)$, introduced by Maejima and Tudor \cite{mt12} as the double Wiener-It\^{o} integral 
\begin{equation*}Z_{\gamma_1,\gamma_2}(t)=\int_{\mathbb{R}^2}^{\prime}\bigg(\int_0^t(s-x_1)_+^{\gamma_1}(s-x_2)_+^{\gamma_2}\,\mathrm{d}s\bigg)\,\mathrm{d}B_{x_1}\,\mathrm{d}B_{x_2},
\end{equation*}
where the prime $\prime$ indicates exclusion of the diagonals $x_1=x_2$ in the stochastic integral, $B_{x}$ is standard Brownian motion and $\gamma_i\in(-1,-1/2)$, $i=1,2$, and $\gamma_1+\gamma_2>-3/2$. 
The Rosenblatt process (see Taqqu \cite{t75}) is the special case $Z_\gamma(t)=Z_{\gamma,\gamma}(t)$, $-3/4<\gamma<-1/2$. It is readily seen that $Z_{\gamma_1,\gamma_2}(t)=_d t^{2+\gamma_1+\gamma_2}Z_{\gamma_1,\gamma_2}(1)$ (see, for example, Gaunt \cite{gaunt vgiii}), and so, for simplicity, we will work with the random variable $Z_{\gamma_1,\gamma_2}(1)$; results for general $t>0$ follow from a rescaling.
For $\phi\in(0,1)$, we define the random variable $Y_\phi$ by
\begin{equation*}Y_\phi=\frac{a_\phi}{\sqrt{2}}(X_1-1)-\frac{b_\phi}{\sqrt{2}}(X_2-1),
\end{equation*}
where $X_1$ and $X_2$ are independent $\chi_{(1)}^2$ random variables and 
\begin{align*}a_\phi=\frac{(2\sqrt{\phi})^{-1}+(\phi+1)^{-1}}{\sqrt{(2\phi)^{-1}+2(\phi+1)^{-2}}}, \quad b_\phi=\frac{(2\sqrt{\phi})^{-1}-(\phi+1)^{-1}}{\sqrt{(2\phi)^{-1}+2(\phi+1)^{-2}}}.
\end{align*}
Suppose $\gamma_1\geq\gamma_2$ and that $\gamma_2=(\gamma_1+1/2)/\phi-1/2$. Then, it was shown by Bai and Taqqu \cite{bt17} that $Z_{\gamma_1,\gamma_2}(1)\rightarrow_d Y_\phi$ as $\gamma_1\rightarrow-1/2$. (Observe that if $\gamma_1\rightarrow-1/2$, then $\gamma_2\rightarrow-1/2$.)

We now observe that by the representation (\ref{gamrep}) we have that $Y_\phi$ is distributed as the product of two correlated normal random variables with $s=(1+\phi)/\sqrt{1+6\phi+\phi^2}$ and $\rho=2\sqrt{\phi}/(\phi+1)$ (here we solved $s(1+\rho)/2=a_\phi/\sqrt{2}$ and $s(1-\rho)/2=b_\phi/\sqrt{2}$). Now, Arras et al.\ \cite{aaps17} showed that, for any $i\geq2$, as $\gamma_1\rightarrow-1/2$,
\begin{equation*}\kappa_i(Z_{\gamma_1,\gamma_2}(1))=\kappa_i(Y_\phi)+O\Big(-\gamma_1-\frac{1}{2}\Big).
\end{equation*}
Inserting this asymptotic relation into (\ref{optvg0}) and (\ref{optvg}) implies that, as $\gamma_1\rightarrow-1/2$, 
\begin{align*}d_{\mathrm{W}}(Z_{\gamma_1,\gamma_2}(1),Y_\phi)&\leq C_\phi\sqrt{-\gamma_1-\frac{1}{2}}, \\
d_{2}(Z_{\gamma_1,\gamma_2}(1),Y_\phi)&\leq C_\phi'\Big(-\gamma_1-\frac{1}{2}\Big).
\end{align*}
where the absolute constants $C_\phi$ and $C_\phi'$ only depend on $\phi$.}

\appendix

\numberwithin{equation}{section}

\setcounter{equation}{0}

\section{The modified Bessel function of the second kind}\label{appa}

In this appendix, we collect some properties of the modified Bessel function of the second kind that are needed in this paper. Unless otherwise stated, these properties can be found in Olver et al.\ \cite{olver}.
The \gr{modified Bessel function of the second kind} $K_\nu(x)$ is defined, for $\nu\in\mathbb{R}$ and $x>0$, by
\[K_\nu(x)=\int_0^\infty \mathrm{e}^{-x\cosh(t)}\cosh(\nu t)\,\mathrm{d}t.
\]
For $x>0$, the function $K_\nu(x)$ is positive for all $\nu\in\mathbb{R}$.
For $\nu=m+1/2$, $m=0,1,2,\ldots$, we have
\begin{equation}\label{special} K_{m+1/2}(x)=\sqrt{\frac{\pi}{2x}}\sum_{j=0}^m\frac{(m+j)!}{(m-j)!j!}(2x)^{-j}\mathrm{e}^{-x}.
\end{equation}
In particular,
\begin{equation}\label{special2}K_{1/2}(x)=\sqrt{\frac{\pi}{2x}}\mathrm{e}^{-x}, \quad\!\! K_{3/2}(x)=\sqrt{\frac{\pi}{2x}}\bigg(1+\frac{1}{x}\bigg)\mathrm{e}^{-x}, \quad\!\!  K_{5/2}(x)=\sqrt{\frac{\pi}{2x}}\bigg(1+\frac{3}{x}+\frac{3}{x^2}\bigg)\mathrm{e}^{-x}.
\end{equation}
The modified Bessel function of the second kind has the following asymptotic behaviour:
\begin{eqnarray}\label{Ktend0}K_{\nu} (x) &\sim& \begin{cases} 2^{|\nu| -1} \Gamma (|\nu|) x^{-|\nu|}, & \quad x \downarrow 0, \: \nu \not= 0, \\
-\log x, & \quad x \downarrow 0, \: \nu = 0, \end{cases} \\
\label{Ktendinfinity} K_{\nu} (x) &\sim& \sqrt{\frac{\pi}{2x}} \mathrm{e}^{-x}, \quad x \rightarrow \infty,\: \nu\in\mathbb{R}.
\end{eqnarray}
A differentiation formula is given by:
\begin{align}
\label{ddbk}\frac{\mathrm{d}}{\mathrm{d}x}\big(x^\nu K_\nu(x)\big)&=-x^{\nu} K_{\nu-1}(x).
\end{align}
The following integral formulas can be found in Gradshetyn and Ryzhik \cite{gradshetyn}. For $a,\alpha>0$ and $\nu>-1/2$,
\begin{equation}\label{besint}\int_0^a t^{\nu}K_{\nu}(\alpha t)\,\mathrm{d}t =\frac{\sqrt{\pi}2^{\nu-1}\Gamma(\nu+1/2)}{\alpha^\nu}a\big(K_{\nu}(\alpha a)\mathbf{L}_{\nu-1}(\alpha a)+K_{\nu-1}(\alpha a)\mathbf{L}_{\nu}(\alpha a)\big),
\end{equation}
\gr{where $\mathbf{L}_\nu(x)$ is a modified Struve function of the first kind (see Chapter 11 of Olver et al.\ \cite{olver}).} For $\mu>|\nu|$ and $\beta<\alpha$,
\begin{align}\int_0^\infty t^{\mu-1}\mathrm{e}^{\beta t}K_\nu(\alpha t)\,\mathrm{d}t&=\frac{\sqrt{\pi}(2\alpha)^\nu}{(\alpha-\beta)^{\mu+\nu}}\frac{\Gamma(\mu+\nu)\Gamma(\mu-\nu)}{\Gamma(\mu+1/2)}\times\nonumber\\
\label{hyint}&\quad\times{}_{2}F_1\bigg(\mu+\nu,\nu+\frac{1}{2};\mu+\frac{1}{2};-\frac{\alpha+\beta}{\alpha-\beta}\bigg).
\end{align}
The following asymptotic approximation follows from a simple rescaling of limiting form (2.13) in Gaunt \cite{gaunt ineq3}. For $\beta<\alpha$ and $\nu>-1/2$,
\begin{equation}\label{kintap}\int_x^\infty\mathrm{e}^{\beta t}t^\nu K_{\nu}(\alpha t)\,\mathrm{d}t\sim\sqrt{\frac{\pi}{2\alpha}}\frac{1}{(\alpha-\beta)}x^{\nu-1/2}\mathrm{e}^{-(\alpha-\beta)x}, \quad x\rightarrow\infty.
\end{equation}
The ratio $K_{\nu-1}(x)/K_\nu(x)$ satisfies the following inequalities.  For $x>0$,
\begin{equation}\label{seg1}\frac{x}{\nu-1/2+\sqrt{(\nu-1/2)^2+x^2}}<\frac{K_{\nu-1}(x)}{K_\nu(x)}<\frac{x}{\nu-1+\sqrt{(\nu-1)^2+x^2}}, \quad \nu>1/2,
\end{equation}
and
\begin{equation}\label{seg3}\frac{K_{\nu-1}(x)}{K_\nu(x)}\leq\frac{x}{\nu-1/2+\sqrt{(\nu-3/2)^2+x^2}}, \quad \nu\geq3/2,
\end{equation}
with equality if and only if $\nu=3/2$. The lower and upper bounds in (\ref{seg1}) were obtained by Segura \cite{segura} and 
Ruiz-Antol\'{i}n and Segura \cite{rs16}, respectively, whilst inequality (\ref{seg3}) can be found in Gaunt and Merkle \cite{gm21}.

\subsection*{Acknowledgements}
The author is supported by a Dame Kathleen Ollerenshaw Research Fellowship.  \gr{I would like to thank the reviewers and associate editor for their excellent suggestions, which helped me to substantially improve the paper.}

\end{document}